\documentclass{amsart}
\usepackage[dvips]{graphicx,psfrag}
\usepackage{amsmath,amsthm,amssymb,IMjournal}
\usepackage{color}
\newcommand{\R}{\mathbb R}
\newcommand{\N}{\mathbb N}

\newcommand{\dd}{\mathrm{d}}
\newcommand{\lam}{\lambda}

\newcommand{\bfx}{\mathbf{x}}
\newcommand{\bfy}{\mathbf{y}}

\newcommand{\bfw}{\mathbf{w}}
\newcommand{\hK}{\widehat{K}}

\newcommand{\hv}{\hat{v}}
\newcommand{\T}{\mathcal{T}}
\newcommand{\I}{\mathcal{I}}

\begin{document}
\newtheorem{theorem}{Theorem}[section]
\newtheorem{definition}[theorem]{Definition}
\newtheorem{lemma}[theorem]{Lemma}
\newtheorem{corollary}[theorem]{Corollary}

\numberwithin{equation}{section}

\title{\textit{A priori}
error estimates for Lagrange interpolation on triangles}
\author{\|Kenta |Kobayashi|, Kunitachi,
        \|Takuya |Tsuchiya|, Matsuyama}

\rec {May 31, 2015}

\abstract
We present the error analysis of Lagrange interpolation on
triangles.  A new \textit{a priori} error estimate is derived in
which the bound is expressed in terms of the diameter and circumradius
of a triangle.  No geometric conditions on triangles are imposed 
in order to get this type of error estimates.
\endabstract



\thanks
The authors are supported by JSPS Grant-in-Aid for Scientific Research
(C) 25400198 and (C) 26400201.
The second author is partially supported by
JSPS Grant-in-Aid for Scientific Research (B) 23340023.
\endthanks

\section{Introduction}

Lagrange interpolation on triangles and the associated error
estimates are important subjects in numerical analysis.
In particular, they are crucial in the error analysis of
finite element methods.  It is well known that we must impose
some geometric condition on the triangles to obtain an error estimation
\cite{BabuskaAziz,  Krizek, Zenisek, Zlamal}.
In the following, we mention some common estimations.

Let $K \subset \R^2$ be an arbitrary triangle with vertices $\bfx_1$,
$\bfx_2$, and $\bfx_3$.  Let $\mathcal{P}_1$ be the set of polynomials with
two variables whose order is at most $1$. For a continuous function
$v \in C^0(\overline{K})$, the Lagrange interpolation
$\I_K^1 v \in \mathcal{P}_1$ of order $1$ is defined by
$v(\bfx_i) = (\I_K^1 v)(\bfx_i)$, $i = 1$, $2$, $3$.  For $K$,
let $h_K$ be its length of the longest edge, and $\rho_K$ be the
diameter of its inscribed circle.

\newcommand{\vvskip}{\vspace{5pt}}
\vvskip
\noindent
\textbf{The minimum angle condition},  Zl\'amal \cite{Zlamal} (1968),
\v{Z}en\'i\v{s}ek \cite{Zenisek} (1969). \\
\textit{
Let $\theta_0$, $(0 < \theta_0 \le \pi/3)$ be a constant. If
any angle $\theta$ of $K$ satisfies $\theta \ge \theta_0$
and $h_K \le 1$, then there exists a constant $C = C(\theta_0)$
independent of $h_K$ such that}
\[
   \|v - \I_K^1 v\|_{1,2,K} \le C h_K |v|_{2,2,K}, \qquad
   \forall v \in H^2(K).
\]

\vvskip
Many textbooks on finite element methods, such as those by
Ciarlet \cite{Ciarlet}, Brenner-Scott \cite{BrennerScott}, and
Ern-Guermond \cite{ErnGuermond}, explain the following theorem.

\vvskip
\noindent
\textbf{Shape-regularity}.
{\it Let $\sigma > 0$ be a constant. If $h_K/\rho_K \le \sigma$ and
$h_K \le 1$, then there exists a constant $C$ that is
independent of $h_K$ such that}
\begin{equation}
 \|v - \I_K^1 v\|_{1,2,K} \le C \frac{h_K^2}{\rho_K}|v|_{2,2,K}
  \le C \sigma h_K |v|_{2,2,K} , \qquad  \forall v \in H^2(K).
  \label{standard-est}
\end{equation}

\vvskip
It is a simple exercise to show that the minimum angle condition is
equivalent to the shape-regularity for triangular elements in $\R^2$.
The maximum of the ratio $h_K/\rho_K$ in a triangulation is called the
\textbf{chunkiness parameter} \cite{BrennerScott}.
The shape-regularity condition is sometimes called the
\textbf{inscribed ball condition} as well. On the conditions
equivalent to the shape-regularity, see \cite{BrandtsKorotovKrizek}.
The minimum angle
condition or shape-regularity, however, are not necessarily needed
to obtain an error estimate.  The following condition is well known.

\vvskip
\noindent
\textbf{The maximum angle condition},
 Babu\v{s}ka-Aziz \cite{BabuskaAziz} (1976). \\
{\it Let $\theta_1$, $(\pi/3 \le \theta_1 < \pi)$ be a constant.  If
any angle $\theta$ of $K$ satisfies $\theta \le \theta_1$ and
$h_K \le 1$, then there exists a constant $C = C(\theta_1)$ 
that is independent of $h_K$ such that}
\[
   \|v - \I_K^1 v\|_{1,2,K} \le C h_K |v|_{2,2,K}, \qquad
   \forall v \in H^2(K).
\]

\vvskip
Later, K\v{r}\'{i}\v{z}ek \cite{Krizek} introduced
the \textit{semiregularity condition}, which is equivalent to
the maximum angle condition (see Section~\ref{concl} (2)).
Let $R_K$ be the circumradius of $K$.

\vvskip
\noindent
\textbf{The semiregularity condition},
K\v{r}\'{i}\v{z}ek \cite{Krizek} (1991). \\
{\it Let $p > 1$ and  $\sigma > 0$ be a constant. If $R_K/h_K \le \sigma$ and
$h_K \le 1$, then there exists a constant $C = C(\sigma)$ 
that is independent of $h_K$ such that}
\[
   \|v - \I_K^1 v\|_{1,p,K} \le C h_K |v|_{2,p,K}, \qquad
   \forall v \in W^{2,p}(K).
\]

\vvskip
Since its discovery, the maximum angle condition has been considered
the most essential condition for error estimates of Lagrange
interpolation on triangular elements.
However, Hannukainen, Korotov and K\v{r}\'i\v{z}ek \cite{HKK} pointed out
that \textit{the maximum angle condition is not necessary for convergence of
the finite element method} by showing simple numerical examples.
Furthermore, the authors recently reported the following error estimation.

\vvskip
\noindent
\textbf{The circumradius condition}, Kobayashi-Tsuchiya
\cite{KobayashiTsuchiya1} (2014). \\
\textit{For an arbitrary triangle $K$ with $R_K \le 1$, there
exists a constant $C_p$ that is independent of $K$
such that the following estimate holds}: 
\begin{equation}
   \|v - \I_K^1 v\|_{1,p,K} \le C_p R_K |v|_{2,p,K}, \qquad
   \forall v \in W^{2,p}(K), \quad 1 \le p \le \infty.
  \label{circumradius-condition1}
\end{equation}

\vvskip
Note that estimate \eqref{circumradius-condition1} follows from
\begin{gather}
    B_p^{1,1}(K) := \sup_{v \in \T_p^{1}(K)}
    \frac{|v|_{1,p,K}}{|v|_{2,p,K}} \le C_p R_K,
  \label{Bp11}
\end{gather}
where the set $\T_p^1(K) \subset W^{2,p}(K)$ is defined by
\[
   \T_p^{1}(K) := \left\{v \in W^{2,p}(K) \Bigm|
     v(\bfx_i) = 0, \; i = 1,2,3 \right\}. 
\]

Suppose that $\{\mathcal{T}_h\}_{h>0}$ is a sequence of triangulations
of a polygonal domain $\Omega\subset\R^2$ such that
\begin{equation}
  \lim_{h\to 0} \max_{K \in \mathcal{T}_h} R_K =0.
  \label{circumradius-condition}
\end{equation}
Let $S_{h}$ be the set of all piecewise linear functions
on $\mathcal{T}_h$, defined by
\[
   S_{h} := \left\{v_h \in H_0^1(\Omega) \cap C(\overline{\Omega})
   \bigm| v_h|_{K} \in \mathcal{P}_1, \forall K \in \mathcal{T}_h  \right\},
\]
and let $u_h \in S_{h}$ be the piecewise linear finite element solution
on the triangulation $\mathcal{T}_h$ of the Poisson problem
\begin{equation*}
  - \Delta u = f \text{ in } \Omega, \qquad
   u = 0 \text{ on }  \partial\Omega
\end{equation*}
for a given $f \in L^2(\Omega)$.  Then, C\'ea's lemma
\cite[Theorem~2.4.1]{Ciarlet} claims that, for the exact solution $u$,
\begin{align*}
    \|u - u_h\|_{1,2,\Omega} \le 
   C \inf_{v_h \in S_{h}}
   | u - v_h |_{1,2,\Omega}
   \le C  | u - \I_{h}^1 u|_{1,2,\Omega}
  \le C  \left(\max_{K \in \mathcal{T}_h} R_K\right) |u|_{2,2,\Omega},
\end{align*}
where $\I_{h}^1 u$ is the global piecewise linear interpolation
of $u$ defined by $\I_{h}^1 u|_K = \I_K^1 u$ for any
$K \in \mathcal{T}_h$.  Hence, if \eqref{circumradius-condition} holds and
$u \in H^2(\Omega)$, the finite element solutions $\{u_h\}$ converge to
$u$ as $h \to 0$.  The condition \eqref{circumradius-condition} is
called the \textbf{circumradius condition} in
\cite{KobayashiTsuchiya1}.

Let $\alpha$, $\beta \in \R$ be such that $1 < \alpha < \beta < 1 + \alpha$.
Consider the triangle $K$ whose vertices are
$(0,0)^T$, $(h,0)^T$, and $(h^\alpha,h^\beta)^T$.
It is straightforward to see $\rho_K = \mathcal{O}(h^\beta)$
and $R_K = \mathcal{O}(h^{1+\alpha - \beta})$.
Hence, if $h \to 0$, the convergence rates which 
\eqref{standard-est} and \eqref{circumradius-condition1}
yield are 
$\mathcal{O}(h^{2-\beta})$ and $\mathcal{O}(h^{1+\alpha - \beta})$,
respectively.  Therefore, \eqref{circumradius-condition1}
gives a better convergence rate than \eqref{standard-est}.  Moreover,
if $\beta \ge 2$, \eqref{standard-est} does not yield convergence while
\eqref{circumradius-condition1} does.
Note that, when $h \to 0$, the maximum angle of $K$ approaches to $\pi$.

From these facts,
we can say that the circumradius $R_K$ of $K$ is more important
than its minimum and maximum angles (or the chunkiness parameter).  It
should also be noted that the circumradius condition is closely related to
the definition of surface area \cite{KobayashiTsuchiya2}.

The aim of this paper is to extend \eqref{circumradius-condition1} to
higher-order Lagrange interpolation and to prove the following theorem.
\begin{theorem} \label{cor6}
Let $K$ be an arbitrary triangle.  Let
$1 \le p \le \infty$, and $k$, $m$ be integers such that $k \ge 1$ and
$0 \le  m \le k$.  Then, for the $k$th-order Lagrange interpolation
 $\I_K^k$ on $K$, the following estimation holds:
\begin{equation}
   |v - \I_K^k v|_{m,p,K} \le C
   \left(\frac{R_K}{h_K}\right)^m h_K^{k+1-m} |v|_{k+1,p,K}
  =  C R_K^m h_K^{k+1-2m} |v|_{k+1,p,K}
  \label{main-estimate}
\end{equation}
for any $v \in W^{k+1,p}(K)$,
where the constant $C$ depends only on $k$, $p$ and is independent of
the geometry of $K$.
\end{theorem}
We here emphasize that 
\textit{no geometric condition on the triangles is imposed}
in Theorem~\ref{cor6}.  Therefore, the estimation \eqref{main-estimate}
is valid even if the maximum angle condition does not hold.

To prove Theorem~\ref{cor6}, we make use of two key observations.  One
of them is that ``\textit{squeezing an isosceles right triangle
perpendicularly does not reduce the approximation property of Lagrange
interpolation,}'' which was first noted by Babu\v{s}ka and Aziz
\cite{BabuskaAziz} for the case $k=1$ and $p=2$. 
This obervation is stated rigorously in Theorem~\ref{squeezedtheorem}.

Note that an arbitrary triangle $K$ can be obtained by
``folding'' or ``unfolding'' an right triangle.  Let $A$
be the $2 \times 2$ matrix that defines the linear transformation of
``folding'' and ``unfolding'' (see \eqref{matrixA}).  Liu and Kikuchi
\cite{LiuKikuchi} pointed out that an error estimation of the linear Lagrange
interpolation $\I_K^1$ is obtained by considering the eigenvalues of
$A^{T}A$.  In Section~\ref{alternate}, we rewrite Liu and
Kikuchi's proofs using Kronecker products of matrices, and one of their
main results \cite[Corollary~1]{LiuKikuchi} is immediately obtained
(Theorem~\ref{liu-kikuchi}).  The other key observation is that the
upper bound in Theorem~\ref{liu-kikuchi} is closely  related to the
circumradius $R_K$ of $K$ (Lemma~\ref{kobayashi}).  Combining
Theorem~\ref{liu-kikuchi} and Lemma~\ref{kobayashi}, an alternative
proof of \eqref{Bp11} is obtained for the case $p=2$
(Corollary~\ref{alternateproof}).

This method is straightforwardly extended to higher-order Lagrange
interpolation in Section~4,  and we obtain the main results of
Theorem~\ref{thm5} that is equivalent to Theorem~\ref{cor6}.

\section{Preliminaries}
\subsection{Notation}
Let $n \ge 1$ be a positive integer and $\R^n$ be the $n$-dimensional
Euclidean space.
Throughout of this paper, $K$ is a triangle in $\R^2$.
We denote the Euclidean norm of $\bfx \in \R^n$ by
$|\bfx|$.  Let $\R^{n*} := \{l:\R^n \to \R : l \text{ is linear}\}$
be the dual space of $\R^n$.  We always regard $\bfx \in \R^n$ as a
column vector and $\mathbf{a} \in \R^{n*}$ as a row vector.
For a matrix $A$ and $\bfx \in \R^n$, $A^T$ and $\bfx^T$
denote their transpositions.
For matrices $A$ and $B$,
$A \otimes B$ denotes their Kronecker product.
For a differentiable function $f$ with $n$ variables,
its gradient $\nabla f = \mathrm{grad} f \in \R^{n*}$ is the row vector
\[
  \nabla f = \nabla_\bfx f := 
  \left(\frac{\partial f}{\partial x_1}, ..., 
    \frac{\partial f}{\partial x_n}\right), \qquad
    \bfx := (x_1, ..., x_n)^T.
\]

Let $\N_{0}$ be the set of nonnegative integers.
For $\delta = (\delta_1,...,\delta_n) \in (\N_{0})^n$,
the multi-index $\partial^\delta$ of partial differentiation 
(in the sense of distribution) is defined by
\[
    \partial^\delta = \partial_\bfx^\delta
    := \frac{\partial^{|\delta|}\ }
   {\partial x_1^{\delta_1}...\partial x_n^{\delta_n}}, \qquad
   |\delta| := \delta_1 + ... + \delta_n.
\]

Let $\Omega \subset \R^n$ be a (bounded) domain.  The usual Lebesgue
space is denoted by $L^p(\Omega)$ for $1 \le p \le \infty$.
For a positive integer $k$, the Sobolev space $W^{k,p}(\Omega)$ is
defined by
$\displaystyle
  W^{k,p}(\Omega) := 
  \left\{v \in L^p(\Omega) \, | \, \partial^\delta v \in L^p(\Omega), \,
   |\delta| \le k\right\}$.
The norm and semi-norm of $W^{k,p}(\Omega)$ are defined,
for $1 \le p < \infty$, by
\begin{gather*}
  |v|_{k,p,\Omega} := 
  \biggl(\sum_{|\delta|=k} |\partial^\delta v|_{0,p,\Omega}^p
   \biggr)^{1/p}, \quad   \|v\|_{k,p,\Omega} := 
  \biggl(\sum_{0 \le m \le k} |v|_{m,p,\Omega}^p \biggr)^{1/p},
\end{gather*}
and $\displaystyle   |v|_{k,\infty,\Omega} := 
  \max_{|\delta|=k} \left\{\mathrm{ess}
   \sup_{\hspace{-5mm}\bfx \in\Omega}|\partial^\delta v(\bfx)|\right\}$,
 $\displaystyle   \|v\|_{k,\infty,\Omega} := 
  \max_{0 \le m  \le k} \left\{|v|_{m,\infty,\Omega}\right\}$.

\subsection{Preliminaries from matrix analysis}
We introduce some facts from the theory of matrix analysis.  For their
proofs, readers are referred to textbooks on matrix analysis such as
\cite{HornJohnson} and \cite{Yamamoto2}.

Let $n \ge 2$ be an integer and $A$ be an $n \times n$ regular matrix. Let
$B:=A^{-1}$. Then, $A^TA$ is symmetric positive-definite and has $n$ positive
eigenvalues.  Let $0 < \mu_m \le \mu_M$ be the minimum and maximum eigenvalues.
Then, we have
\[
  \mu_m |\bfx|^2 \le  |A\bfx|^2 \le \mu_M |\bfx|^2, \quad
  \mu_M^{-1} |\bfx|^2 \le  |B^T\bfx|^2  \le \mu_m^{-1}
   |\bfx|^2,  \quad
    \forall \bfx \in \R^n.
\]
Then, the minimum and maximum eigenvalues of 
$(A^TA)\otimes(A^TA) = (A\otimes A)^T(A\otimes A)$ are
$0 < \mu_m^2 \le \mu_M^2$.  Hence, for any $\bfw\in \R^{n^2}$, we have
\begin{align*}
   \mu_m^2 |\bfw|^2 \le |(A \otimes A)\bfw|^2
   & \le \mu_M^2 |\bfw|^2,  \quad
  \mu_M^{-2} |\bfw|^2 \le |(B\otimes B)^T\bfw|^2
   \le \mu_m^{-2} |\bfw|^2.
\end{align*}

The above facts can be straightforwardly extended to the case of
the higher-order Kronecker product $A\otimes ... \otimes A$.
For $A\otimes ... \otimes A$, $B\otimes ... \otimes B$
(the $k$th Kronecker products), and we have, for $\bfw \in \R^{n^k}$, 
\begin{align*}
   \mu_m^k |\bfw|^2  \le |(A \otimes ... \otimes  A)\bfw|^2
   \le \mu_M^k |\bfw|^2, \;
  \mu_M^{-k} |\bfw|^2 \le |(B\otimes ... \otimes B)^T\bfw|^2
   \le \mu_m^{-k} |\bfw|^2.
\end{align*}

\subsection{The affine transformation defined by a regular matrix}
Let $A$ be an $n \times n$ matrix with det$A > 0$.
We consider the affine transformation $\varphi(\bfx)$ defined by
$\bfy = \varphi(\bfx) := A \bfx + \mathbf{b}$ for
   $\bfx = (x_1, ..., x_n)^T$, $\bfy = (y_1, ..., y_n)^T$ with
   $\mathbf{b} \in \R^n$.
Suppose that a reference region $\widehat\Omega\subset \R^n$ is
transformed to a domain $\Omega$ by $\varphi$;
$\Omega := \varphi(\widehat\Omega)$.   Then, a function $v(\bfy)$
defined on $\Omega$ is pulled-back to the function $\hv(\bfx)$
on $\widehat\Omega$ as    $\hv(\bfx) := v(\varphi(\bfx)) = v(\bfy)$.
Then, we have $\nabla_\bfx \hv = (\nabla_\bfy v) A$,
$\nabla_\bfy v = (\nabla_\bfx \hv) B$, and
$|\nabla_\bfy v|^2 = |(\nabla_\bfx \hv) B|^2
    = (\nabla_\bfx \hv) BB^T (\nabla_\bfx \hv)^T$.

The Kronecker product $\nabla\otimes\nabla$ of the gradient $\nabla$ is
defined by
\[
  \nabla\otimes\nabla := \left(
   \frac{\partial\ }{\partial x_1}\nabla, ...,
   \frac{\partial\ }{\partial x_n}\nabla\right)
 = \left(
   \frac{\partial^2\ }{\partial x_1^2},
   \frac{\partial^2\ }{\partial x_1\partial x_2}, ...,
   \frac{\partial^2\ }{\partial x_{n-1}\partial x_n},
   \frac{\partial^2\ }{\partial x_n^2}\right).
\]
We regard $\nabla\otimes\nabla$ to be a row vector.
From this definition, it follows that
\[
   \sum_{|\delta| = 2} (\partial^\delta v)^2 = \sum_{i,j=1}^n 
  \left(\frac{\partial^2 v}{\partial x_i\partial x_j}\right)^2
  = |(\nabla\otimes\nabla) v|^2
\]
and  $(\nabla_\bfx\otimes\nabla_\bfx) \hv =
   \left((\nabla_\bfy\otimes\nabla_\bfy) v\right) (A\otimes A)$, 
   $(\nabla_\bfy\otimes\nabla_\bfy) v =
   \left((\nabla_\bfx\otimes\nabla_\bfx) \hv \right) (B\otimes B)$.
Suppose that the minimum and maximum eigenvalues of $BB^T$
are $0 < \lam_m \le \lam_M$.  Then, we have
$\lam_m |\nabla_\bfx \hv|^2  \le
   |\nabla_\bfy v|^2 \le \lam_M |\nabla_\bfx \hv|^2$ and
\begin{align*}
  \sum_{|\delta| = 2} (\partial_\bfy v)^2 & = 
  |(\nabla_\bfy\otimes\nabla_\bfy) v|^2 \\
 & = \left((\nabla_\bfx\otimes\nabla_\bfx) \hv\right) (B\otimes B)
       (B\otimes B)^T \left((\nabla_\bfx\otimes\nabla_\bfx) \hv\right)^T  \\
 & = \left((\nabla_\bfx\otimes\nabla_\bfx) \hv\right) (BB^T\otimes BB^T)
       \left((\nabla_\bfx\otimes\nabla_\bfx) \hv\right)^T, \\
  \lam_m^2\sum_{|\delta| = 2} (\partial_\bfx^\delta \hv)^2 
 & \le \sum_{|\delta| = 2} (\partial_\bfy^\delta v)^2 
 \le \lam_M^2 \sum_{|\delta| = 2} (\partial_\bfx^\delta \hv)^2.
\end{align*}
The above inequalities can be easily extended to higher-order
derivatives giving the following inequalities:
\begin{align}
  \lam_m^k\sum_{|\delta| = k} (\partial_\bfx^\delta \hv)^2 
 \le \sum_{|\delta| = k} (\partial_\bfy^\delta v)^2 
 \le \lam_M^k \sum_{|\delta| = k} (\partial_\bfx^\delta \hv)^2, \qquad k \ge 1.
  \label{generalm}
\end{align}

\subsection{Useful inequalities}
For $N$ positive real numbers $U_1, ..., U_N$, the following
inequalities hold:
\begin{gather}
   \sum_{k=1}^N U_k^p \le N^{\tau(p)}
   \left(\sum_{k=1}^N U_k^2\right)^{p/2}, \quad
   \tau(p) := \begin{cases}
         1-p/2, & 1 \le p \le 2 \\
         0,     & 2 \le p < \infty
     \end{cases},  \label{tau}\\
    \left(\sum_{k=1}^N U_k^2\right)^{p/2} \le N^{\gamma(p)}
   \sum_{k=1}^N U_k^p, \quad
   \gamma(p) := \begin{cases}
         0, &  1 \le p \le 2 \\
         p/2 - 1,     & 2 \le p < \infty
     \end{cases}.
    \label{gamma}
\end{gather}

\subsection{The Sobolev imbedding theorems}
If $1 < p < \infty$, Sobolev's Imbedding Theorem and Morry's inequality
imply that
\begin{gather*}
   W^{2,p}(K) \subset C^{1,1-2/p}(\overline{K}), \quad  p > 2, \\
   H^{2}(K) \subset W^{1,q}(K) \subset C^{0,1-2/q}(\overline{K}),
     \quad \forall q > 2, \\ 
   W^{2,p}(K) \subset W^{1,2p/(2-p)}(K) \subset C^{0,2(p-1)/p}(\overline{K}),
    \quad  1 < p < 2.
\end{gather*}
For proofs of the Sobolev imbedding theorems, see
\cite{AdamsFournier} and \cite{Brezis}.
For the case $p=1$, we still have the continuous imbedding
   $W^{2,1}(K) \subset C^0(\overline{K})$.
For the proof of the critical imbedding , see
\cite[Theorem~4.12]{AdamsFournier} and \cite[Lemma~4.3.4]{BrennerScott}.

\subsection{Lagrange interpolation on triangles and their estimations}
Let $K$ be a triangle with vertices
$\bfx_i$, $i=1,2,3$, and $(\lam_1,\lam_2,\lam_3)$ be its
barycentric coordinates with respect to $\bfx_i$.
By definition, we have $0 \le \lam_i \le 1$,
$\lam_1+\lam_2+\lam_3=1$.  For a positive integer $k \ge 1$,
the set $\Sigma^k(K)$ of points on $K$ is defined by
\begin{equation}
   \Sigma^k(K) := \left\{\left(
   \frac{a_1}{k},\frac{a_2}{k},\frac{a_3}{k}\right) \in K \Bigm|
    a_i \in \N_0, \; 0 \le a_i \le k,\; a_1+a_2+a_3 = k \right\}.
  \label{Sigmak}
\end{equation}
\begin{figure}[thb]
\begin{center}
  \includegraphics[width=10cm]{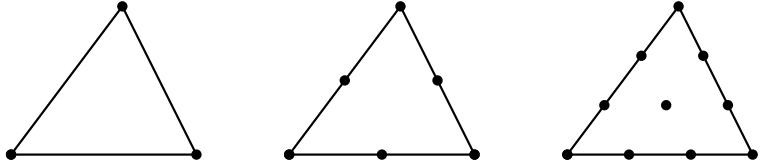} 
\caption{The set $\Sigma^k(K)$, $k=1$, $k=2$, $k=3$.}
 \label{gridpoints}
\end{center}
\end{figure}

For a triangle $K$, a positive integer $k$, and
$1 \le p \le \infty$, we define the subset
$\T_p^{k}(K) \subset W^{k+1,p}(K)$ by
\begin{gather}
   \T_p^{k}(K) := \left\{v \in W^{k+1,p}(K) \Bigm|
     v(\bfx) = 0, \; \forall \bfx \in \Sigma^k(K) \right\}.
  \label{TpkK}
\end{gather}
Let $\mathcal{P}_k$ be the set of polynomials with two variables whose
degree is at most $k$.  For a continuous function
$v \in C(\overline{K})$, the $k$th-order Lagrange interpolation
$\I_K^k v \in \mathcal{P}_k$ of $v$ is defined by
$v(\bfx) = (\I_K^kv)(\bfx)$ for any $\bfx \in \Sigma^k(K)$.
From this definition, it is clear that
     $v - \I_K^k v \in \T_p^{k}(K)$ for any
     $v \in W^{k+1,p}(K)$.

For an integer $m$ such that $0 \le m \le k$, $B_p^{m,k}(K)$
is defined by
\begin{gather*}
    B_p^{m,k}(K) := \sup_{v \in \T_p^{k}(K)}
    \frac{|v|_{m,p,K}}{|v|_{k+1,p,K}}.
\end{gather*}
Note that we have
\begin{gather}
 B_p^{m,k}(K) = \inf\left\{C ; 
   |v - \I_K^k v|_{m,p,K} \le C |v|_{k+1,p,K}, \
   \forall v \in W^{k+1,p}(K)\right\}.
 \label{equivalent}
\end{gather}

For an error estimate of Lagrange interpolation, standard textbooks
such as \cite{Ciarlet} and \cite{BrennerScott} explain the following
theorem.  Recall that $\rho_K$ is the diameter of its inscribed circle of $K$.
\begin{theorem}
Let $1 \le p \le \infty$, and $k \ge 1$ be an integer.
Let $\sigma > 0$ be a positive constant.
Then, for a triangle $K$ that satisfies $h_K/\rho_K \le \sigma$, 
the following estimate holds:
\begin{equation}
  |v - \I_K^k v|_{m,p,K} \le C h_K^{k+1-m} |v|_{k+1,p,K},
  \quad \forall v \in W^{k+1,p}(K),
  \label{standard}
\end{equation}
where $m=0,1,...,k$, and the constant $C$ depends on $k$, $p$,
and $\sigma$.
\end{theorem}

Jamet presented an improved estimation, which does not require the
shape-regularity condition \cite[Th\'eor\`eme~3.1]{Jamet}.
\begin{theorem}[Jamet]\label{Jamettheorem}
Let $1 \le p \le \infty$. Let $m \ge 0$, $k \ge 1$ be integers such that
$k+1-m > 2/p$ $(1 < p \le \infty)$ or $k-m \ge 1$ $(p=1)$.
\footnote{Note that in \cite[Th\'{e}or\`{e}me~3.1]{Jamet}
the case $p=1$ is not mentioned explicitly but clearly holds for
triangles.}
Then, the following estimate holds:
\begin{equation}
  |v - \I_K^k v|_{m,p,K} \le C \frac{h_K^{k+1-m}}{(\cos(\theta_K/2))^m}
  |v|_{k+1,p,K},   \quad \forall v \in W^{k+1,p}(K),
  \label{jamet}
\end{equation}
where $\theta_K \ge \pi/3$ is the maximum angle of $K$, and
$C$ depends only on $k$, $p$.
\end{theorem}

Note that, if $m=k \ge 1$ and $1 \le p \le 2$, estimate \eqref{jamet}
cannot be applied.  As will be noted in Section~\ref{concl} (2), 
Theorem~\ref{cor6} includes Theorem~\ref{Jamettheorem} as a special
case.  

Let $K_\alpha$ be the right triangle with vertices $(0,0)^T$, $(1,0)^T$,
and $(0,\alpha)^T$ $(0 < \alpha \le 1)$, that is obtained by
squeezing $\hK$.  As is stated in Section~1, squeezing
a right triangle perpendicularly does not deteriorate approximation
property of Lagrange interpolation.  We have the following theorem:

\begin{theorem}\label{squeezedtheorem}
There exists a constant $C_{k,p}$ that depends only on
$k$ and $p$ $(1 \le p \le \infty)$ and is independent of 
$\alpha$ $(0 < \alpha \le 1)$ such that
\begin{equation}
   B_p^{m,k}(K_\alpha) := \sup_{v \in \T_p^{k}(K_\alpha)}
  \frac{|v|_{m,p,K_\alpha}}{|v|_{k+1,p,K_\alpha}}
   \le C_{k,p}, \qquad m = 0,1,...,k.
  \label{extension}
\end{equation}
\end{theorem}


Note that Theorem~\ref{squeezedtheorem} is not a totally new result.
For the case $m=k=1$ and $p=2$, \eqref{extension} was proved by
Babu\v{s}ka and Aziz in \cite{BabuskaAziz}. 
Kobayashi and Tsuchiya \cite{KobayashiTsuchiya1} proved \eqref{extension}
with $m=k=1$ and any $p$ $(1 \le p \le \infty)$.
For the case $k \ge 1$ with $p=2$ and $m=0$, $1$,
\eqref{extension} was proved by Shenk \cite{Shenk}.  By \eqref{jamet},
estimate \eqref{extension} holds if
$k+1-m > 2/p$ $(1 < p \le \infty)$ or $k - m \ge 1$ $(p=1)$. 
Hence, it seems that \eqref{extension} with $k=m \ge 2$ and
$1 \le p \le 2$ has not yet been proved.
A proof of Theorem~\ref{squeezedtheorem} by the Babu\v{s}ka-Aziz
type technique will be given in \cite{KobayashiTsuchiya3}.

\section{Liu and Kikuchi's method} \label{alternate}
In this section, we give an alternative proof of \eqref{Bp11} for the case
$p=2$ using Liu and Kikuchi's method. To this end, we rewrite their proof
using the Kronecker product of matrices.

\begin{figure}[thb]
\begin{center}
  \psfrag{a}[][]{$e_1$}
  \psfrag{b}[][]{$e_3$}
  \psfrag{c}[][]{$e_2$}
  \psfrag{t}[][]{$\theta$}
  \psfrag{x1}[][]{$\bfx_1$}
  \psfrag{x2}[][]{$\bfx_2$}
  \psfrag{x3}[][]{$\bfx_3$}
  \psfrag{K}[][]{$K$}
  \includegraphics[width=6cm]{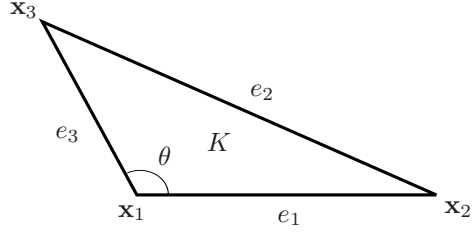}
 \caption{The triangle under consideration.
 The vertices are
 $\bfx_1=(0,0)^T$, $\bfx_2=(1,0)^T$, and $\bfx_3=(\alpha s,\alpha t)^T$,
 where $s^2 + t^2 = 1$, $t > 0$, and $0 < \alpha \le 1$.
 We assume that $|e_1| = 1 \le |e_2|$.}
 \label{LiuKikuchiTriangle}
\end{center}
\end{figure}

For $s$, $t$, and $\alpha$ with $s^2 + t^2 = 1$, $t > 0$,
$0 < \alpha \le 1$, we consider the vector
$(\alpha s,\alpha t)^T \in \R^2$.  Let $K$ be the triangle
with vertices $\bfx_1:=(0,0)^T$,
$\bfx_2:=(1,0)^T$, and $\bfx_3:=(\alpha s,\alpha t)^T$.
Let $e_1$, $e_2$, $e_3$ be the three edges of $K$, as depicted in
Figure~\ref{LiuKikuchiTriangle}.
Without loss of generality, we assume
that $e_2$ is the longest edge of $K$.
Let $\theta$ be the angle between $e_1$ and $e_3$.
Then, $s = \cos\theta$, $t = \sin\theta$, and the assumption
that $e_2$ is the longest yields
\begin{equation}
   s= \cos\theta \le \frac{\alpha}{2} \le \frac{1}{2}, \qquad
    \frac{\pi}{3} \le \theta < \pi.
  \label{conditionK}
\end{equation}
Note that an arbitrary triangle in $\R^2$ can be transformed to $K$ by a
sequence of scaling, translation, rotation, and mirror imaging.

We define the $2 \times 2$ matrices as 
\begin{align}
  A := \begin{pmatrix}
	   1 &  s \\ 0 &  t 
	 \end{pmatrix}, \qquad
    B := A^{-1} = \begin{pmatrix}
	   1 & - st^{-1} \\ 0 & t^{-1}
	 \end{pmatrix}.
  \label{matrixA}
\end{align}
Then, $K_\alpha$ can be transformed to $K$ by the transformation
$\bfy = A\bfx$.  Moreover, $\T_p^k(K)$ is pulled-back to
$\T_p^k(K_\alpha)$ as
 $\T_p^k(K) \ni v \mapsto \hv := v\circ A \in\T_p^k(K_\alpha)$.
A simple computation yields that $A^TA$ has eigenvalues $1 \pm |s|$, and
$BB^T$ has eigenvalues $(1\mp|s|)/t^2$.  It follows from
\eqref{generalm} that
$\frac{1-|s|}{t^2} |\nabla_\bfx \hv|^2 \le
   |\nabla_\bfy v|^2 \le \frac{1+|s|}{t^2} |\nabla_\bfx \hv|^2$ and
\begin{gather}
  \frac{(1-|s|)^2}{t^4}\sum_{|\delta|=2}
   \left(\partial_\bfx^\delta \hv \right)^2
 \le  \sum_{|\delta| = 2} \left(\partial_\bfy^\delta v\right)^2  \le
  \frac{(1+|s|)^2}{t^4}\sum_{|\delta|=2}
  \left(\partial_\bfx^\delta \hv \right)^2.
  \label{k=2}
\end{gather}
Furthermore, because the determinant of $A$ is $t$, we have
\begin{gather*}
   |v|_{1,2,K}^2 \le \frac{1+|s|}{t} |\hv|_{1,2,K_\alpha}^2, \qquad
  \frac{(1-|s|)^2}{t^3} |\hv|_{2,2,K_\alpha}^2 \le |v|_{2,2,K}^2,
   \\
   \frac{|v|_{1,2,K}^2}{|v|_{2,2,K}^2} \le
   \frac{t^2(1+|s|)|\hv|_{1,2,K_\alpha}^2}
     {(1-|s|)^2|\hv|_{2,2,K_\alpha}^2}
  = \frac{(1+|s|)^2|\hv|_{1,2,K_\alpha}^2}
     {(1-|s|)|\hv|_{2,2,K_\alpha}^2}.
\end{gather*}
Combining this estimate and \eqref{extension} with $m=k=1$ and $p=2$,
we obtain the following theorem \cite[Corollary~1]{LiuKikuchi}:
\begin{theorem}[Liu-Kikuchi]  \label{liu-kikuchi}
For $0 < \alpha \le 1$, we have the following estimate:
\[
 B_2^{1,1}(K) \le \frac{1+|s|}{\sqrt{1-|s|}}B_2^{1,1}(K_\alpha)
   \le  \frac{2 C_{1,2}}{\sqrt{1-|s|}}.
\]
\end{theorem}

The following is the key lemma.

\begin{lemma} \label{kobayashi}
Let $R_K$ be the circumradius of $K$.  For the triangle $K$ considered
in this section, the following inequality holds:
\[
    \frac{1}{\sqrt{1-|s|}} \le 2 \sqrt{2}\, R_K.
\]
\end{lemma}
\proof
 Recall from \eqref{conditionK} that $s = \cos\theta$,
$t = \sin\theta$, and $\pi/3 \le \theta < \pi$.
A straightforward computation implies that
\begin{equation*}
   \sqrt{1+|s|} \le \sqrt{2}\sqrt{1 + \alpha^2 - 2 \alpha s}, \qquad
    0 < \forall \alpha \le 1, \; -1 < \forall s \le \frac{\alpha}{2}.
\end{equation*}
From the cosine and sine laws,
we have $|e_2|^2 = 1 + \alpha^2 - 2\alpha s = 4R_K^2 t^2$.
Therefore, we obtain 
\begin{align*}
   \frac{1}{\sqrt{1-|s|}} = \frac{\sqrt{1+|s|}}{t}
    \le \frac{\sqrt{2}}{t} \,\sqrt{1 + \alpha^2 - 2 \alpha s} 
    =  \frac{\sqrt{2}}{t} \,\sqrt{4R_K^2 t^2}
    = 2 \sqrt{2}\, R_K. 
\end{align*}
\endproof

Combining Theorem~\ref{liu-kikuchi} and Lemma~\ref{kobayashi}, we
have obtained an alternative proof of \eqref{Bp11} for the
triangle depicted in Figure~\ref{LiuKikuchiTriangle} with $p=2$.

\begin{corollary}\label{alternateproof}
Let $K$ be the triangle depicted in Figure~\ref{LiuKikuchiTriangle}.
Then, we have
\[
    B_2^{1,1}(K) := \sup_{v\in \T_2^{1}(K)} 
    \frac{|v|_{1,2,K}}{|v|_{2,2,K}} \le 4 \sqrt{2}\, C_{1,2} R_K.
\]
\end{corollary}

\section{Main results and their proofs}
The method explained so far can be immediately extended to higher-order
Lagrange interpolation.  Inequality \eqref{k=2} is extended to
the case of arbitrary $k$ as follows:
\[
  \frac{(1-|s|)^k}{t^{2k}}\sum_{|\delta|=k} (\partial_\bfx^\delta \hat{v})^2
 \le   \sum_{|\delta|=k} (\partial_\bfy^\delta v)^2
 \le \frac{(1+|s|)^k}{t^{2k}} \sum_{|\delta|=k} (\partial_\bfx^\delta \hat{v})^2.
\]

Let $1 \le p < \infty$.
Then, inequalities \eqref{tau}, \eqref{gamma} yield
{\allowdisplaybreaks
\begin{align*}
  |v|_{m,p,K}^p & = \int_K \sum_{|\delta|=m}
            |\partial_\bfy^\delta v(\bfy)|^p \dd \bfy
    \le 2^{m \tau(p)} \int_K \left(
 \sum_{|\delta|=m} |\partial_\bfy^\delta v(\bfy)|^2 \right)^{p/2}
    \hspace{-3.5mm} \dd \bfy \\
   & \le 2^{m \tau(p)} 
   \left(\frac{1+|s|}{t^2}\right)^{mp/2}
      \int_K \left(
 \sum_{|\delta|=m} |\partial_\bfx^\delta \hv(\bfx)|^2 \right)^{p/2}
    \hspace{-3.5mm} \dd \bfy \\
   & = 2^{m \tau(p)} 
   \left(\frac{1+|s|}{t^2}\right)^{mp/2}t
 \int_{K_\alpha} \left(
  \sum_{|\delta|=m} |\partial_\bfx^\delta \hv(\bfx)|^2 \right)^{p/2}
   \hspace{-3.5mm}      \dd \bfx \\
   & \le 2^{m (\tau(p)+\gamma(p))} 
   \left(\frac{1+|s|}{t^2}\right)^{mp/2}t
      \int_{K_\alpha} 
      \sum_{|\delta|=m} |\partial_\bfx^\delta \hv(\bfx)|^p  \dd \bfx \\
    & = 2^{m (\tau(p)+\gamma(p))} 
   \left(\frac{1+|s|}{t^2}\right)^{mp/2}t |\hv|_{m,p,K_\alpha}^p
\end{align*}
}
and
{\allowdisplaybreaks
\begin{align*}
  |v|_{k+1,p,K}^p & = \int_K \sum_{|\delta|=k+1}
            |\partial_\bfy^\delta v(\bfy)|^p \dd \bfy \\
   & \ge 2^{- (k+1) \gamma(p)} \int_K \left(
      \sum_{|\delta|=k+1} |\partial_\bfy^\delta v(\bfy)|^2 \right)^{p/2}
    \hspace{-3.5mm} \dd \bfy \\
   & \ge 2^{ -(k+1) \gamma(p)} 
   \left(\frac{1-|s|}{t^2}\right)^{(k+1)p/2}
      \int_K \left(
         \sum_{|\delta|=k+1} |\partial_\bfx^\delta \hv(\bfx)|^2 \right)^{p/2}
     \hspace{-3.5mm} \dd \bfy \\
   & = 2^{ -(k+1) \gamma(p)} 
   \left(\frac{1-|s|}{t^2}\right)^{(k+1)p/2}t
      \int_{K_\alpha} \left(
     \sum_{|\delta|=k+1} |\partial_\bfx^\delta \hv(\bfx)|^2 \right)^{p/2}
    \hspace{-3.5mm} \dd \bfx \\
   & \ge 2^{-(k+1) (\tau(p)+\gamma(p))} 
   \left(\frac{1-|s|}{t^2}\right)^{(k+1)p/2}t
      \int_{K_\alpha} 
      \sum_{|\delta|=k+1} |\partial_\bfx^\delta \hv(\bfx)|^p  \dd \bfx \\
    & = 2^{-(k+1) (\tau(p)+\gamma(p))} 
   \left(\frac{1-|s|}{t^2}\right)^{(k+1)p/2}t |\hv|_{k+1,p,K_\alpha}^p.
\end{align*}
}
The two inequalities and Theorem~\ref{squeezedtheorem}, 
Lemma~\ref{kobayashi} imply
\begin{align*}
    \frac{|v|_{m,p,K}^p}{|v|_{k+1,p,K}^p} & \le
  \tilde{c}_{k,m,p}^p \frac{t^{p(k+1-m)}(1+|s|)^{mp/2}|\hv|_{m,p,K_\alpha}^p}
     {(1-|s|)^{(k+1)p/2}|\hv|_{k+1,p,K_\alpha}^p} \\
 &  = \tilde{c}_{k,m,p}^p \frac{(1+|s|)^{(k+1+m)p/2}|\hv|_{m,p,K_\alpha}^p}
     {t^{pm}|\hv|_{k+1,p,K_\alpha}^p}, \\
   \frac{|v|_{m,p,K}}{|v|_{k+1,p,K}} & \le \tilde{c}_{k,m,p} 
   \frac{(1+|s|)^{(k+1+m)/2}|\hv|_{m,p,K_\alpha}}{t^m
 |\hv|_{k+1,p,K_\alpha}}  \le c_{k,p} C_{k,p} R_K^m,
\end{align*}
where $\tilde{c}_{k,m,p} := 2^{(k+1+m)(\tau(p)+\gamma(p))/p}$ and the
constant $c_{k,p}$ depends only on $k$, $p$.
If $p = \infty$, the
same estimation is obtained by letting $p \to \infty$ in the above
inequalities.
Thus, denoting $c_{k,p} C_{k,p}$ by $C_{k,p}$, the
following theorem has been proved.

\begin{theorem}\label{thm4}
Let $K$ be the triangle depicted in Figure~\ref{LiuKikuchiTriangle}.
Then, the estimate
\[
    B_p^{m,k}(K) := \sup_{v\in \T_p^{k}(K)} 
    \frac{|v|_{m,p,K}}{|v|_{k+1,p,K}} \le \, C_{k,p} R_K^m,
      \qquad \forall p, \ 1 \le p \le \infty
\]
holds, where $R_K$ is the circumradius of $K$ and the
constant $C_{k,p}$ depends only on $k$ and $p$.
\end{theorem}

Now, let $K$ be an arbitrary triangle.
Theorem~\ref{thm4} and Corollary~\ref{alternateproof} can be extended to
$K$.  The similar transformation $G_Y$ for a positive
$Y \in \R$ is defined by $G_Y:\R^2 \to \R^2$, $G_Y({\bf x}) := Y {\bf x}$.
Let $K_1$ be defined by $K_1 = G_Y(K)$.
A function $u \in W^{k,p}(K)$ on $K$ is pulled-back to
$v(\bfx) := u(G_Y^{-1}(\bfx)) = u(G_{1/Y}(\bfx))$ on
$K_1$.  Then, for a nonnegative integer $k$ and any
$p$ $(1 \le p \le \infty)$, we have
\[
  |v|_{k,p,K_1} = Y^{2/p-k}|u|_{k,p,K}, \qquad
  \forall u \in W^{p,k}(K).
\]

Let $h_K \ge h_2 \ge h_1$ be the lengths of the three edges of $K$.
Suppose that the second longest edge of $K$ is parallel to
the $x$- or $y$-axis.  Then, by a translation,
a mirror imaging, and $G_{1/h_2}$, $K$ can be transformed to the
triangle $\widetilde{K}$ depicted in Figure~\ref{LiuKikuchiTriangle}.
Hence, we may apply Theorem~\ref{thm4} to $\widetilde{K}$, and obtain
\begin{gather*}
    \sup_{u\in \T_p^{k}(K)} 
    \frac{h_2^{m}|u|_{m,2,K}}{h_2^{k+1}|u|_{k+1,2,K}} =
    \sup_{v\in \T_p^{k}(\widetilde{K})} 
    \frac{|v|_{m,p,\widetilde{K}}}
         {|v|_{k+1,p,\widetilde{K}}} \le \, C_{k,p}
     R_{\widetilde{K}}^m
\end{gather*}
and
\begin{gather*}
    \sup_{u\in \T_p^{k}(K)} 
    \frac{|u|_{m,p,K}}{|u|_{k+1,p,K}}
      \le \, C_{k,p} R_{\widetilde{K}}^m h_2^{k+1-m}
      \le \, C_{k,p} R_{K}^m h_K^{k+1-2m}.
\end{gather*}
Here, we use the fact that $R_{\widetilde{K}}h_2 = R_K$ and
$h_K/2 < h_2 \le h_K$.  The constant $C_{k,p}$ can be
modified up to a constant multiple.
Note that if $p \neq 2$, the Sobolev norms are modified by a rotation.
Therefore, we have shown the following theorem, which is equivalent to
Theorem~\ref{cor6} because of \eqref{equivalent}.

\begin{theorem}\label{thm5}
Let $K$ be an arbitrary triangle.  Let $R_K$ be its
circumradius and $h_K$ be the length of its longest edge.  Let
$1 \le p \le \infty$, and $m$, $k$ be integers such that
$0 \le  m \le k$.   Then, there exists a positive constant
$C$ that depends only on $k$, $p$ such that
the following estimation holds:
\[
    B_p^{m,k}(K) := \sup_{u\in \T_p^{k}(K)} 
    \frac{|u|_{m,p,K}}{|u|_{k+1,p,K}} \le \, C
   \left(\frac{R_K}{h_K}\right)^m h_K^{k+1-m}
   = C R_K^m h_K^{k+1-2m}.
\]
\end{theorem}

\subsection{Concluding remarks}\label{concl}
Here, we compare the newly obtained estimate \eqref{main-estimate}
with known results such as \eqref{standard}, \eqref{jamet}, and
\eqref{circumradius-condition1}.

(1) For an error analysis of the finite element method, the cases
$m=0, 1$ are the most important.  In these cases, the estimates obtained
from \eqref{main-estimate} can be written, for any $v \in W^{k+1,p}(K)$, as
\begin{gather*}
   |v - \I_K^k v|_{1,p,K} \le C R_K h_K^{k-1} |v|_{k+1,p,K}, \qquad
    |v - \I_K^k v|_{0,p,K} \le C h_K^{k+1} |v|_{k+1,p,K}.
\end{gather*}
They are extensions of \eqref{circumradius-condition1}.
Recall that the constant $C$ is independent of the geometry
of $K$.

(2) Recall that $h_1 \le h_2 \le h_K$ are the lengths of the
three edges of $K$.  Let $\theta_K$ be the maximum angle of $K$ and
$S_K$ be the area of $K$.
Then, from the formulas $S_K = \frac{1}{2}h_1h_2\sin\theta_K$ and
$R_K = h_1h_2h_K/(4S_K)$, we have
\[
    \frac{R_K}{h_K}  = \frac{1}{2\sin\theta_K}, \qquad
     \frac{\pi}{3} \le \theta_K < \pi.
\]
Thus, it is clear that
\textit{the boundedness of $R_K/h_K$, which is the semiregularity
of $K$ defined by K\v{r}\'i\v{z}ek, is equivalent to
the maximum angle condition $\theta_K \le \theta_1 < \pi$} with a
fixed constant $\theta_1$.  If this is the case, the estimate
from \eqref{main-estimate} becomes
\[
   \left|v - \I_K^k v \right|_{m,p,K} \le 
   \frac{C}{(2 \sin \theta_1)^m} h_K^{k+1-m} |v|_{k+1,p,K}, \quad
    \forall v \in W^{k+1,p}(K)
\]
for $m = 0,1,...,k$, which is an extention of Jamet's result of
\eqref{jamet}.

{\small
{\em Kenta Kobayashi}, 
          Graduate School of Commerce and Management, 
       Hitotsubashi University, Kunitachi, Japan
       e-mail:  \texttt{kenta.k@r.hit-u.ac.jp}. \\
{\em Takuya Tsuchiya},
       Graduate School of Science and Engineering,
          Ehime University, Matsuyama, Japan
          e-mail: \texttt{tsuchiya@math.sci.ehime-u.ac.jp}.
}
\end{document}